\documentclass[12pt]{article}
\usepackage[a4paper,left=2.5cm,right=2.5cm,top=2.5cm,bottom=2.5cm]{geometry}

\usepackage{enumitem}

\usepackage{graphicx}
\usepackage{amsfonts}
\usepackage{amsmath}
\usepackage{amssymb}
\usepackage[ansinew]{inputenc}
\usepackage{latexsym}
\usepackage{tabularx}
\usepackage{makeidx}
\usepackage{color}
\allowdisplaybreaks
\newtheorem{op}{Open Problem}

\newtheorem{con}{Conjecture}

\newtheorem{theorem}{Theorem}[section]

\newcommand{\NN}{{\mathbb N}}

\newcommand{\ZZ}{{\mathbb Z}}
\newcommand{\RR}{{\mathbb R}}

\newcommand{\bsz}{\boldsymbol{z}}
\newcommand{\bsx}{\boldsymbol{x}}

\newcommand{\bsy}{\boldsymbol{y}}

\newcommand{\bsf}{\boldsymbol{f}}
\newcommand{\bsg}{\boldsymbol{g}}

\newcommand{\bsS}{{\bf S}}

\newcommand{\bsT}{{\bf T}}

\newcommand{\bsalpha}{\boldsymbol{\alpha}}

\newcommand{\To}{\rightarrow}

\begin{document}

\title{\scshape{Discrepancy estimates for sequences: new results and open problems}}
\author{Gerhard Larcher\thanks{is partially supported by the Austrian Science Fund (FWF), Project P21943.}     }
\date{}
\maketitle

\abstract{ In this paper we give an overview of recent results on (upper and lower) discrepancy estimates for (concrete) sequences in the unit-cube. 
In particular we also give an overview of discrepancy estimates for certain classes of hybrid sequences. Here by a hybrid sequence we understand an $(s+t)$-dimensional sequence which is a combination of an $s$-dimensional sequence of a certain type (e.g. Kronecker-, Niederreiter-, Halton-, $\ldots$ type) and a $t$-dimensional sequence of another type. The analysis of the discrepancy of hybrid sequences (and of their components) is a rather current and vivid branch of research.
We give a collection of some challenging open problems on this topic.}
\vspace{1em}


\section{Introduction}

Let $(\bsz_n)_{n \ge 0}$ be a sequence in the $d$-dimensional unit-cube $\left[0,1\right)^d $. Then the {\it discrepancy} of the first $N$ points of the sequence is defined by
$$D_N = \sup_{B \subseteq \left[0,1\right)^d }\left|\frac{A_N(B)}{N}-\lambda (B)\right|,$$
where 
$$A_N(B) := \;\#\; \{n : 0 \le n < N,  \bsz_n \in B\}, $$
$\lambda$ is the $d$-dimensional volume and the supremum is taken over all axis-parallel subintervals $B \subseteq \left[0,1)^d.\right.$
The sequence $(z_n)_{n \ge 0}$ is called uniformly distributed if $\lim\limits_{N \rightarrow\infty} D_N = 0.$ 
If the supremum is restricted to all $B$ with the left-lower corner in the origin, then we speak of the star-discrepancy $D_N^*$.
It is the most well-known conjecture in the theory of irregularities of distribution, that for every sequence $(z_n)_{n \ge 0}$ in $\left[0,1)^d \right.$ we have 
$$ D_N \ge c_d  \frac{(\log N )^d }{N}$$
for a constant $c_d > 0$ and for infinitely many $N$. This result was proven for dimension $d=1$ by Schmidt \cite{Schmidt2}.
Hence sequences whose discrepancy satisfies $ D_N = O \left( \frac{(\log N )^d }{N}\right)$ are called ``low-discrepancy sequences''. \newline
Note that recent investigations of Bilyk, Lacey et al., see for example \cite{Bilyk} or \cite{Bilyk et}, have led some people to conjecture that $\frac{(\log N )^{\frac{d+1}{2}} }{N}$ instead of $\frac{(\log N )^d }{N}$ is the best possible order for the discrepancy of sequences in $\left[0,1\right)^d$.  
At the moment the best known general lower bound for the discrepancy of sequences in $\left[0,1\right)^d$ for $d \ge 2$ is 
$$D_N \ge c_d  \frac{(\log N )^{\frac{d}{2}+\epsilon (d)} }{N}$$
for infinitely many $N$, with some small $\epsilon (d) > 0.$ For more details on this topic see \cite{Beck2} or \cite{BilykL}. 

There are three groups of (almost) low-discrepancy sequences which are of main interest for applications in quasi-Monte Carlo methods. (Here by a quasi-Monte Carlo method we mean simulation in the setting of Monte Carlo methods, but using deterministic, i.e., quasi-random (usual low-discrepancy) point sets instead of pseudo-random point sets.) Indeed, these are (until now) the only known types of sequences containing concrete examples of (almost) low-discrepancy sequences. 
These are Kronecker sequences, Halton sequences and digital $(\bsT,s)$-sequences in the sense of Niederreiter. 

The most classical type are the {\em Kronecker sequences}. A Kronecker sequence is of the form
$$ \bsz_n = \left(\left\{n \bsalpha\right\}\right)_{n \ge 0} = \left(\left(\left\{n \alpha_1\right\}, \ldots, \left\{n \alpha_d\right\}\right)\right)_{n \ge 0}$$
for some $\bsalpha =\left(\alpha_1,\ldots,\alpha_d\right) \in \left[0,1)^d\right.$.  
The sequence is uniformly distributed in $\left[0,1)^d \right.$  iff $1,\alpha_1, \ldots , \alpha_d$ are linearly independent over $\ZZ$.

(Good) {\it lattice point sets} are -- in some sense -- finite versions of Kronecker sequences. They are of the following form:
$$ \bsz_n = \left(\left\{ n \frac{a_1}{N}\right\}, \ldots, \left\{n \frac{a_d}{N}\right\}\right)_{n=0, 1 \ldots, N-1}$$
for some given $N \in \NN$ and $a_1, \ldots, a_d \in \ZZ$.

The second type of sequences are {\em Halton sequences} which are defined as follows: 
Let $b_1, \ldots , b_d \ge 2$ be pairwise coprime integers, then the Halton sequence $\left(\bsz_n\right)_{n \ge 0}$ with $\bsz_n=\left(z_n^{(1)},\ldots,z_n^{(d)}\right)$, in bases $b_1, \ldots , b_d$ is given by 
$$z_n^{(j)} := \psi_{b_j}(n)$$
where $$\psi_{b_j}(n) := \sum_{i=0}^{\infty}n_i b_j^{-i-1}$$ for $n \in \NN_0$ with base $b_j$ representation
$$n = \sum_{i=0}^{\infty} n_i b_j^i \quad, 0 \le n_i \le b_j \, .$$
It is easy to show that every Halton sequence is a low-discrepancy sequence. 

The concept of {\em digital-sequences} in a base $q$ (see for example \cite{Dick}, \cite{Lar}, \cite{Nied1}, or \cite{Nied2}) was introduced by Niederreiter and it contains also earlier special examples of sequences of the same type considered by Sobol' and by Faure.

Let $q$ be a prime number and let $\ZZ_q$ be the finite field of order $q$. We identify $\ZZ_q$ with the set $\{0,1,\ldots, q-1\}$ equipped with arithmetic operations modulo $q$. 
Let $d \in \NN$. Let $C_1, \ldots, C_d \in \ZZ_q^{\NN \times \NN}$ be $\NN \times \NN$ matrices over $\ZZ_q$. Let $n \in \NN_0$, with $q$-adic expansion $n = n_0 + n_1q + n_2q^2 + \ldots$ and set
$$ \vec{n} = (n_0, n_1, n_2, \ldots)^\top \in (\ZZ_q^\NN)^\top.$$
Then define
$$ \vec{x}_{n,j} = C_j \vec{n} \quad \mbox{for} \quad j=1, \ldots, d, $$
where all arithmetic operations are taken modulo $q$. Let $\vec{x}_{n,j} = (x_{n,j,1}, x_{n,j,2}, \ldots)^\top$ and define
$$ x_{n,j} = x_{n,j,1}q^{-1} + x_{n,j,2}q^{-2} + \cdots \, .$$
Then the $n$th point $\bsx_n$ of the sequence $\bsS (C_1, \ldots, C_d)$ is given by $\bsx_n = (x_{n,1}, \ldots, x_{n,d})$. A sequence $\bsS (C_1, \ldots, C_d)$ constructed this way is called a \textit{digital sequence (over $\ZZ_q$)} with generating matrices $C_1, \ldots, C_d$.
Under certain conditions on the generating matrices it can be shown that the star discrepancy of $\bsS (C_1, \ldots, C_d)$ is of order of magnitude $\frac{(\log N)^d}{N}$ in $N$. For more information we refer to \cite{Dick}, \cite{Nied1}, \cite{Nied2} and the references therein.

The finite versions of digital sequences are digital $(t,m,s)$-nets (over $\ZZ_q$). These are point-sets $\bsx_n = (x_{n,1}, \ldots, x_{n,d})$ with $n=0,1,\ldots, N-1$, where $N=q^m$, and which are defined in the same way as digital sequences  but now with $C_1, \ldots, C_d \in \ZZ_q^{m \times m}$.

In the following we will give some recent results on the discrepancy of these point sequences and point sets, and of combinations (hybrids) of them, and we will also give a collection of challenging open problems in this context.

\section{Metrical and average type discrepancy estimates for digital point sets and sequences and for good-lattice point-sets}

In \cite{Beck} Beck has given upper and lower metrical bounds for the discrepancy of the $d$-dimensional Kronecker sequence. He showed:

\begin{theorem}(Beck, 1994)\label{Beck1}
For almost all $\bsalpha$ in $[0,1)^d$ for the discrepancy $D_N$ of the $d$-dimensional Kronecker sequence we have
$$ D_N = O \left(\frac{(\log N)^d (\log \log N)^{1+\epsilon}}{N}\right) $$
for every $\epsilon > 0$.
\end{theorem}

\begin{theorem}(Beck, 1994)\label{Beck2}
For almost all $\bsalpha$ in $[0,1)^d$ for the discrepancy $D_N$ of the $d$-dimensional Kronecker sequence we have 
$$ D_N \ge C (\alpha, d) \frac{(\log N)^d  \log \log N}{N} $$
for infinitely many $N$, where $C (\alpha, d) > 0 \cdots $.
\end{theorem} 

It is quite interesting that until now, for dimensions $d \ge 2$ no \textbf{concrete} choice of $\alpha_1, \ldots, \alpha_d$ is known, such that the upper discrepancy estimate in Theorem \ref{Beck1} of Beck is satisfied for this \textbf{concrete} sequence.
The discrepancy of a Kronecker sequence generated by $\bsalpha = (\alpha_1, \ldots, \alpha_d)$ heavily depends on how good $\alpha_1, \ldots, \alpha_d$ can be simultaneously approximated by rationals.
For example if $\alpha_1, \ldots, \alpha_d$ are algebraic numbers such that $1, \alpha_1, \ldots, \alpha_d$ are linearly independent over $\ZZ$, then by the Theorem of Thue-Siegel-Roth-Schmidt we have $D_N = O \left(\frac{1}{N^{1-\epsilon}}\right)$ for every $\epsilon > 0$. Further discrepancy estimates for the Kronecker sequence in dependence on diophantine approximation properties of $\alpha_1, \ldots, \alpha_d$ can be found in \cite{KuNied2006} or in \cite[Theorem 2]{Nied2012}.

For the proofs of Theorems \ref{Beck1} and \ref{Beck2} Beck uses a Poisson summation formula for the discrepancy function and some results from probabilistic diophantine approximation.
Especially, for example, he uses a result of Schmidt from \cite{Schmidt} which gives for almost all $\bsalpha = (\alpha_1, \ldots, \alpha_d) \in [0,1)^d$ a rather exact formula for
\begin{align*}
N(h):= \# & \left\{ (q_1, \ldots, q_d) \in \ZZ^d \, \Big|  \, \left| q_i \right| \le h \quad \mbox{for} \quad i=1, \ldots d \right. \quad \mbox{and}\\
  \quad &  \left\{ q_1 \alpha_1 + \cdots + q_d \alpha_d\right\} < \phi (q_1, \ldots, q_d)\Big\}.
\end{align*}
Indeed, for suitable $\phi$ we have 
$$ N(h) = \underset{\vline \, q_i \, \vline \, \le h}{\underset{q_1, \ldots, q_d}{\sum}}\phi(q_1, \ldots, q_d) + R, $$
with a certain ``small'' error-term R.

In \cite{Lar-Pill} an analogous result to Theorem \ref{Beck1} was given for digital sequences.

\begin{theorem} (Larcher, 1998)
Let $d \in \NN$, let $q$ be a prime number and let $\epsilon > 0$. Then for $\mu_d$-almost all $d$-tuples $(C_1, \ldots, C_d) \in (\ZZ_q^{\NN \times \NN})^d$ of generating matrices the digital sequence generated by $(C_1, \ldots, C_d)$ has discrepancy satisfying
$$ D_N = O \left(\frac{(\log N)^d (\log \log N)^{2+\epsilon}}{N}\right)$$
for all $\epsilon > 0$.
(Here $\mu_d$ is a probability measure on $(\ZZ_q^{\NN \times \NN})^d$ defined in a quite natural way, see \cite{Lar-Pill}.)
\end{theorem}

For the proof of this result one had to combine some counting arguments with results from metrical non-archimedean diophantine approximation.
An important subclass of the class of digital sequences is the class of digital Kronecker sequences. These sequences build a ``non-archimedean analog'' to classical Kronecker sequences. They have been introduced by Niederreiter \cite[Section 4]{Nied2}, and further investigated by Larcher and Niederreiter \cite{Lar2}.

Let $\ZZ_q[x]$ be the set of all polynomials over $\ZZ_q$ and let $\ZZ_q((x^{-1}))$ be the field of formal Laurent series $g$ with $g=0$ or
$$ g = \sum_{k=\omega}^{\infty} a_k x^{-k} \quad \mbox{with} \quad a_k \in \ZZ_q \quad \mbox{and} \quad \omega \in \ZZ \quad \mbox{with} \quad a_\omega \neq 0. $$
$\ZZ_q ((x^{-1}))$ contains the field of rational functions over $\ZZ_q$ as a subfield. The discrete exponential evaluation $\nu$ of $g$ is defined by $\nu(g):=-\omega$ and $\nu(0):=-\infty$. Furthermore, we define the ``fractional part'' of $g$ by
$$ \{g\}:= \sum_{k=\max (1, \omega)}^{\infty} a_k x^{-k}.$$
We associate a nonnegative integer $n$ with $q$-adic expansion $n=n_0 + n_1 q + \cdots + n_rq^r$ with the polynomial $n(x) = n_0 + n_1 x + \cdots + n_r x^r$ in $\ZZ_q [x]$ and vice versa.

For every $d$-tuple $\bsf = (f_1, \ldots, f_d)$ of elements of $\ZZ_q((x^{-1}))$ we define the sequence $S(\bsf) = (\bsx_n)_{n \ge 0}$ by
$$ \bsx_n = (\{n(x)f_1(x)\}_{x=q}, \ldots, \{n(x) f_d (x)\}_{x=q}) \quad \mbox{for} \quad n \in \NN_0. $$
In analogy to classical Kronecker sequences it has been shown in \cite{Lar2} that a digital Kronecker sequence $S(\bsf)$ is uniformly distributed in $[0,1)^d$ if and only if $1,f_1, \ldots, f_d$ are linearly independent over $\ZZ_q [x]$.
By $\mu$ we denote the normalized Haar-measure on $\ZZ_q((x^{-1}))$ and by $\tilde{\mu}_d$ the $d$-fold product measure on $(\ZZ_q((x^{-1})))^d$.

In \cite{Lar1995} Larcher proved the following metrical upper bound on the star discrepancy of digital Kronecker sequences.

\begin{theorem}(Larcher, 1995). 
Let $d \in \NN$, let $q$ be a prime number and let $\epsilon > 0$. For $\tilde{\mu}_d$-almost all $\bsf \in (\ZZ_q((x^{-1})))^d$ the digital Kronecker sequence $S (\bsf)$ has star discrepancy satisfying
$$ D_N (S(\bsf)) = O \left( \frac{(\log N)^d (\log \log N)^{2+\epsilon}}{N}\right). $$ 
\end{theorem}

Quite recently Larcher and Pillichshammer were able to give corresponding metrical lower bounds for the discrepancy of digital sequences and of digital Kronecker sequences (see \cite{Lar-Pill2}, \cite{Lar-Pill3}).

\begin{theorem} (Larcher and Pillichshammer, 2013). \label{Lar-Pill1}
Let $d \in \NN$ and let $q$ be a prime number. Then for $\mu_d$-almost all $d$-tuples $(C_1, \ldots, C_d) \in (\ZZ_q^{\NN \times \NN})^d$ of generating matrices the digital sequence $S (C_1, \ldots, C_d)$ over $\ZZ_q$ has discrepancy satisfying
$$ D_N (S (C_1, \ldots, C_d)) \ge c(q, d)\frac{(\log N)^d \log \log N}{N} \quad \mbox{for infinitely many} \, N \in \NN $$
with some $c (q, d) > 0$ not depending on $N$.
\end{theorem}

\begin{theorem} (Larcher and Pillichshammer, 2013).\label{Lar-Pill2}
Let $d \in \NN$ and let $q$ be a prime number. For $\tilde{\mu}_d$-almost all $\bsf \in (\ZZ_q((x^{-1})))^d$ the digital Kronecker sequence $S (\bsf)$ has star discrepancy satisfying
$$ D_N (S (\bsf)) \ge c(q, d)\frac{(\log N)^d \log \log N}{N} \quad \mbox{for infinitely many} \, N \in \NN $$
with some $c(q,d) > 0$ not depending on $N$.
\end{theorem}

For the proofs of Theorems \ref{Lar-Pill1} and \ref{Lar-Pill2} we used, in analogy to the method of Beck, a Poisson summation formula based on Walsh-functions for the discrepancy-function of these digital sequences and again certain results on non-archimedean diophantine approximation.

For the finite versions of these point sequences, namely for the good lattice point sets (as discrete versions of Kronecker sequences) and for digital $(t,m,s)$-nets (as discrete versions of digital sequences), as well as for digital nets generated by rational functions over finite fields until now just \textit{upper} ``metrical'' (i.e, average type), bounds are known.

The last mentioned class of point set is a finite analogon to the digital Kronecker sequences.
Digital nets generated by rational functions over finite fields are defined as follows:
Let $f \in \ZZ_q [x]$ with $\deg (f) = t \ge 1$, and $g_1, \ldots, g_d \in \ZZ_q [x]$ with $\gcd (g_i, f) = 1$ for $i=1, \ldots, d$ and $\deg (g_i) < t$.

Then we consider point sets $(\bsx_n)_{n \ge 0}$ in $[0,1)^d$ of the form
$$ \bsx_n := \left(\left. \left\{\frac{n(x)g_1(x)}{f(x)}\right\} \right|_{x=q}, \ldots, \left.\left\{\frac{n(x) g_d (x)}{f(x)}\right\}\right|_{x=q}\right).$$
The mentioned upper bounds for the discrepancy of these three classes of finite point sets were proven in \cite{Bykov}, \cite{Lar1993} and \cite{Lar1998}.

\begin{theorem}(Bykovskii, 2012).\label{Bykovskii}
For all $\chi<1$ there is a constant $C_\chi > 0$ such that for all $N$ the discrepancy $D_N$ of the lattice point set
$$ \left(\left\{n \frac{a_1}{N}\right\}, \ldots, \left\{n \frac{a_d}{N}\right\}\right)_{n=0, \ldots, N-1} $$
satisfies 
$$D_N \le C_\chi \frac{(\log N)^{d-1} \log \log N}{N}$$ 
for at least $\chi N^d$ choices for $(a_1, \ldots, a_d)$ with $0 \le a_1, \ldots, a_d < N$.
\end{theorem}

In dimension $d=2$ this result already had been shown by Larcher in \cite{Lar1986}. In dimension $d=2$ the result moreover has a narrow connection to the conjecture of Zaremba on continued fractions. One version of this conjecure is the following:

\begin{con}(Zaremba).
\textit{There is an absolute constant $A$ such that for all $N \in \NN$ there is an $a \in \NN$ relatively prime to $N$ such that all continued fraction coefficients of $\frac{a}{N}$ are less than $A$.}
\end{con}

A weaker version of this conjecture is the following conjecture (which is said to be stated by Moser for the first time):

\begin{con}(Moser)
\textit{There is an absolute constant $B$ such that for all $N \in \NN$ there is an $a \in \NN$ relatively prime to $N$ such that the sum of all continued fraction coefficients of $\frac{a}{N}$ is less than $B \log N$.}
\end{con}

Of course the conjecture of Moser is true if the conjecture of Zaremba holds.
From the correctness of Moser's conjecture it follows that for all $N$ there exists an $a$ relatively prime to $N$ such that for the discrepancy $D_N$ of the $2$-dimensional lattice point set $$\left(\left\{n \frac{1}{N}\right\}, \left\{n \frac{a}{N}\right\}\right)_{n=0, \ldots, N-1}$$
we have
$$ D_N \le C(B) \frac{\log N}{N} ,$$
which is an improvement of Theorem \ref{Bykovskii} for dimension $2$, if we see this just as an existence result.
Concerning the conjecture of Zaremba there has been important progress quite recently by the paper \cite{BouKont} of Bourgain and Kontorovich who show, that the conjecture of Zaremba holds at least for a set of integers $N$ with density $1$. See in this connection also the papers of Frolenkov and Kontorovich \cite{Frolenkov} and \cite{Kontorovich}.

For digital nets and for digital nets generated by rational functions over finite fields we have the following upper average type estimates (see \cite{Lar1993} and \cite{Lar1998}).

\begin{theorem} (Larcher, 1998)\label{Lar2}
For given integers $d \ge 1, m \ge 2$ and prime base $q$ we have: for all $\delta$ with $0 < \delta < 1$, the number of $d$-tuples $C = (C_1, \ldots, C_d)$ of $m \times m$-matrices, providing a digital net over $\ZZ_q$ with discrepancy $D_N$ satisfying
$$ D_N \le \frac{1}{\delta} B(d, q)\frac{(\log N)^{d-1} \log \log N}{N} + O \left(\frac{(\log N)^{d-1}}{N}\right), $$
is at least
$$ (1-\delta) \# M_d (m). $$
(Here $B(d,q)$ is a constant depending only on $d$ and $q$, whereas the O-constant also depends on $\delta$, and $M_d(m)$ is the set of all $d$-tuples of $m \times m$ matrices over $\ZZ_q$.)
\end{theorem}

\begin{theorem}(Larcher, 1993)\label{Lar3}
For every $t \in \NN$ there are $g_1, \ldots, g_d \in \ZZ_q[x], g_1 =1, \gcd (g_i, x) = 1, i=1, \ldots, d$, such that for the discrepancy $D_N$ of the point set
$$ \bsx_n:= \left( \left. \left\{\frac{n(x)g_1(x)}{x^t}\right\} \right|_{x=q} , \ldots, \left. \left\{\frac{n(x)g_d(x)}{x^t}\right\} \right|_{x=q}\right), $$
$$ \quad \quad n=0, \ldots, q^t - 1 =: N-1, $$
we have 
$$ D_N^* < c (d, q) \frac{(\log N)^{d-1} (\log \log N)}{N}, $$
with a constant $c (d, q)$ depending only on $d$ and $q$. 
\end{theorem}

Indeed by the proof of Theorem \ref{Lar3} in \cite{Lar1993} it is easy to see, that this existence result even also can be stated as an average-type result like Theorem \ref{Lar2}.
An analogous result for such point sets, but with denominator  $f \in \ZZ_q[x]$ of degree $t$, with  $\gcd (f,x)= 1$ instead of denominator $x^t$ was shown by Kritzer and Pillichshammer in \cite{KriPill}.

\begin{op}\label{latpois}
\textit{Show that for every $\chi < 1$ there is a constant $C'_{\chi} > 0$ such that for all $N$ for the discrepancy $D_N$ of the lattice point set
$$ \left(\left\{n \frac{a_1}{N}\right\}, \ldots, \left\{n \frac{a_d}{N}\right\}\right)_{n = 0, \ldots, N-1} $$
we have
$$ D_N \ge C'_\chi \frac{(\log N)^{d-1} \log \log N}{N} $$
for at least $\chi N^d$ choices for $(a_1, \ldots, a_d)$ with $0 \le a_1, \ldots, a_d < N$.}
\end{op}

\begin{op}\label{constant}
\textit{Show that for given integers $d \ge 1, m \ge 2$ and prime base $q$ there exists a constant $\tilde{B} (d, q)> 0$ such that: for all $\delta$ with $0 < \delta < 1$, the number of $d$-tuples $C = (C_1, \ldots, C_d)$ of $m \times m$-matrices, providing a digital net over $\ZZ_q$ with discrepancy $D_N$ satisfying
$$ D_N \ge \delta \tilde{B}(d, q) \frac{(\log N)^{d-1} \log \log N}{N}$$
for infinitely many $N$, is at least $(1-\delta) \# M_d(m)$, where $M_d(m)$ is the set of all $d$-tuples of $m \times m$ matrices over $\ZZ_q$.}
\end{op}

\begin{op}\label{ratfunc}
\textit{Give an average-type lower bound (in the style as stated in the Open Problem \ref{constant}) for the discrepancy of digital nets generated by rational functions over finite fields.}
\end{op}

For the proof of Problem \ref{latpois} (and analogously for the proof of Problems \ref{constant} and \ref{ratfunc}) we think that it is possible to use the basic method from the proof of Theorem \ref{Beck2} (respectively from Theorems \ref{Lar-Pill1} and \ref{Lar-Pill2}) together with a certain average type result on discrete diophantine approximation. For example a result of following type, which would be a discrete version of the above mentioned result of Schmidt from \cite{Schmidt} would be quite helpful:

\begin{op}\label{asymprep}
\textit{Find (for suitable $\phi$) an asymptotic representation with a ``small'' error-term $R_\chi$ for  
$$N (h; a_1, \ldots, a_d):=$$ 
$$\# \left\{ -h \le n_1, \ldots, n_d \le h \left| \right. \left\{n_1 \frac{a_1}{N} + n_2 \frac{a_2}{N} + \ldots + n_d \frac{a_d}{N} \right\} < \phi (n_1, \ldots, n_d)\right\}$$
of the following form:
For all $\chi < 1$ and all $h \in \NN$ we have 
$$ N(h; a_1, \ldots, a_d) = \sum_{-h \le n_1, \ldots, n_d < h}  \phi (n_1, \ldots, n_d) + R_\chi(h)$$
for at least $\chi N^d$ choices of $(a_1, \ldots, a_d)$ with $0 \le a_1, \ldots, a_d < N$.}
\end{op}

\noindent For the proof of Open Problems \ref{constant} and \ref{ratfunc} we again probably will need a non-archimedean version of Problem \ref{asymprep}.

It probably should be a quite challenging task to extend the investigations of the above type to Halton-Niederreiter sequences. By a Halton-Niederreiter sequence we understand any $d$-dimensional sequence $(\bsz_n)_{n \ge 0}$ which is obtained by combining $l$ sequences $(\bsz_n^{(1)})_{n \ge 0}, \cdots, (\bsz_n^{(l)})_{n \ge 0}$ in dimensions $d_1, \cdots, d_l$ with $d_1 + \cdots + d_l = d$, i.e.,
$$ \bsz_n = (\bsz_n^{(1)}, \bsz_n^{(2)}, \ldots, \bsz_n^{(l)}), $$
where $(\bsz_n^{(i)})_{n \ge 0}$ is a digital sequence in base $q_i$, and $\gcd (q_1, \cdots, q_l)=1$. 

The Halton sequence is a special case of a Halton-Niederreiter sequence.
It was shown in \cite{HofKriLarPill} that a Halton-Niederreiter sequence is uniformly distributed if and only if each of its $l$ components is uniformly distributed.
The discrepancy of such sequences was investigated in \cite{HofLar} and in \cite{Hofer}. It turns out that, on the one hand, the class of Halton-Niederreiter sequences contains low-discrepancy sequences, but on the other hand in general Halton-Niederreiter sequences are not low-discrepancy sequences even when all of their components are of low-discrepancy. A simple example of such a uniformly distributed Halton-Niederreiter sequence, generated by two low-discrepancy digital sequences, which itself is not of low-discrepancy, is the two-dimensional sequence generated by the matrix 
\[
C_1 =
\begin{pmatrix}
1 & 1 & 1 & 1 & 1 & 1 & \ldots\\
0 & 1 & 0 & 0 & 0 & 0 & \ldots\\
0 & 0 & 1 & 0 & 0 & 0 & \ldots\\
0 & 0 & 0 & 1 & 0 & 0 & \ldots\\
0 & 0 & 0 & 0 & 1 & 0 & \ldots\\
\vdots & \vdots & \vdots & \vdots & \vdots & \vdots \\
\end{pmatrix}
\in \ZZ_3^{\NN \times \NN}
\]
in base 3 and the unit matrix in base 2. 

It would be of interest to determine the order of discrepancy of almost all Halton-Niederreiter sequences, i.e., to solve the following:

\begin{op}\label{copbas}
\textit{For given $d_1, \ldots, d_l \in \NN$ with $d_1 + \ldots + d_l = d$ and given coprime bases $q_1, \ldots, q_l$ determine $g (N)$ as small as possible and $f(N)$ as large as possible such that we have: for the discrepancy $D_N$ of almost all Halton-Niederreiter sequences with $d_i$-dimensional components in base $q_i; i=1, \ldots, l$
we have
$$ D_N = \Omega (f(N)) \quad \mbox{and} \quad D_N = O (g(N)).$$
}\end{op}

We think that, to solve this problem, the techniques developed by Hellekalek in \cite{Hell} based on certain function systems should be quite helpful (see also \cite{Hell2} and \cite{HellNied}.)

In the analysis of Halton-Niederreiter sequences carried out until now, it turned out that this analysis is much easier if the generating matrices of the components all are ``finite-row-matrices''. 
That means: each row of each generating matrix has only finitely many entries different from zero.
Moreover, it seems that the metrical investigation of discrepancy of Halton-Niederreiter sequences will lead to a smaller order of discrepancy if we restrict ourselves to considering ``finite row digital Halton-Niederreiter sequences'' than in the general case. Of course, we first have to consider what a suitable measure for these finite row sequences is.

\begin{op}\label{copbascomp}
\textit{For given $d_1, \ldots, d_l$ with $d_1 + \ldots + d_l = d$ and given coprime bases $q_1, \ldots, q_l$ determine $g (N)$ as small as possible and $f(N)$ as large as possible such that we have: for the discrepancy $D_N$ of almost all \textbf{finite-row} Halton-Niederreiter sequences (with respect to a suitable measure) with $d_i$-dimensional components in base $q_i; i=1, \ldots, l$
we have
$$D_N = \Omega (f(N)) \quad \mbox{and} \quad D_N = O (g(N)).$$
}\end{op}

\section{Discrepancy estimates for and applications of hybrid sequences}\label{discest}
As already mentioned above there are three groups of ((almost) low-discrepancy) sequences which are of main interest. Indeed, these are (until now) the only known types of sequences containing concrete examples of (almost) low-discrepancy sequences. These are Kronecker sequences, Halton sequences and digital sequences in the sense of Niederreiter.
(In some sense we could consider the classes of sequences introduced by Levin \cite{Levin} as a fourth class of low-discrepancy sequences.)

In the last years there grew considerable interest in the distribution of "hybrids" of such sequences, and also of such sequences together with pseudo-random-number-sequences (as well from a theoretical point of view as from the point of view of applications).
A hybrid sequence is defined as follows:
take an $s$-dimensional sequence $(\bsx_n)_{n \ge 0}$ of a certain type and a $t$-dimensional sequence $(\bsy_n)_{n \ge 0}$ of another type and combine them to an $(s+t)$-dimensional {\em hybrid sequence} $$(\bsz_n)_{n \ge 0}  : = ((\bsx_n,\bsy_n))_{n \ge 0}.$$

For the application of these hybrid sequences in QMC methods see for example \cite{Keller}, or \cite{Spanier}. 
A possible, quite natural application of hybrid sequences in financial risk management is described at the end of this section.

Hybrids of the two classical types of sequences, namely of Halton sequences and of Kronecker sequences (we call them Halton-Kronecker sequences) were first studied by Niederreiter in \cite{Nied2009}, \cite{Nied2010} or \cite{Nied2012}.
From a metrical point of view the discrepancy of Halton-Kronecker sequences was studied in \cite{HofLar2} and in \cite{Lar2013}, where the following was shown:

\begin{theorem}(Larcher, 2013).\label{Lar2013}
For every Halton sequence in $[0,1)^s$ and almost every $\boldsymbol{\alpha} \in \mathbb{R}^t$ for the discrepancy of the ($s+t$)-dimensional Halton-Kronecker sequence we have $D_N = O \left( \frac{(\log N)^{s+t + \epsilon}}{N} \right)$ for every $\epsilon > 0$.
\end{theorem}

The proof of this Theorem again heavily depends on the techniques developed by Beck in \cite{Beck}.

So for almost all $\boldsymbol{\alpha} \in \mathbb{R}^t$ a Halton-Kronecker sequence is an (almost) low-discrepancy sequence. However, until now, we do not know any \textbf{{concrete}} example of an (almost) low-discrepancy Halton-Kronecker sequence. Even in the most simple case where $s=t=1$ we do not have any such concrete example. So for example it would be of great interest to study the discrepancy $D_N$ of $ (\bsz_n)_{n \ge 0} = (x_n, \{n \sqrt{2}\})_{n \ge 0} $
where $(x_n)_{n \ge 0}$ is the one-dimensional Halton sequence in base $2$, i.e., the van der Corput sequence in base 2.

\begin{op}\label{vanCorpseq}
\textit{Study the discrepancy of concrete Halton-Kronecker sequences, as for example $(x_n, \{n \sqrt{2}\})_{n \ge 0}$ with $(x_n)_{n \ge 0}$ the van der Corput sequence.}
\end{op}

To handle Open Problem \ref{vanCorpseq} it is necessary to study the growth of the largest coefficients $A_K$ in the continued fraction expansion of $2^K \sqrt{2}$ for $ K=1,2,\ldots $.
So as a prework for solving Problem \ref{vanCorpseq} it would be helpful to give an answer to the following question.

\begin{op}\label{sharboun}
\textit{Let $A_K$ denote the largest continued fraction coefficient of $2^K \sqrt{2}$ for $K=1,2,\ldots$. Give sharp bounds for the growth behaviour of 
$$ B_L := \underset{K \le L}{\max} \; A_K. $$}
\end{op}

Indeed for a solution of Problem \ref{vanCorpseq} even more subtle investigations on the continued fraction coefficients of $2^K \sqrt{2}$ will be necessary.

In \cite{Kritzer} Kritzer has shown an analog to Theorem \ref{Lar2013} for \textbf{finite} hybrid point sets. He showed an ``average-type'' upper bound for the combination of a Hammersley point set with a good lattice point set.
We recall, that a Hammersley point set is a set in $[0,1)^d$ of the form
$$ (\bsx_n)_{n=0}^{N-1}= \left(\frac{n}{N}, \bsy_n\right)_{n=0}^{N-1}$$
where $(\bsy_n)_{n=0}^{N-1}$ are the first $N$ elements of a $(d-1)$-dimensional Halton sequence.

\begin{theorem}(Kritzer, 2012).\label{Kritzer1}
Let $b_1, \ldots, b_s$ be $s$ distinct prime numbers and let $N$ be a prime that is different from $b_1, \ldots, b_s$. Let $H_N := (\bsx_n)_{n=0}^{N-1}$ be the $(s+1)$-dimensional Hammersley point set to the bases $b_1, \ldots, b_s$. Then there exists a generating vector $\bsg \in \{1, \ldots, N-1\}^t$ such that the point set
$$ S_N = (\bsz_n)_{n=0}^{N-1} = ((\bsx_n, \bsy_n))_{n=0}^{N-1} $$
in $[0,1)^{s+t+1}$ with $\bsy_n = \{\frac{n\bsg}{N}\},$ for $0 \le n \le N - 1$, satisfies
$$ D_N (S_N) = O \left( \frac{(\log N)^{s+t+1}}{N}\right), $$
with an implied constant independent of $N$.
\end{theorem}

Indeed, the result of Kritzer is not a discrete analog to Theorem \ref{Lar2013}, but to the result given in \cite{HofLar}, which is a predecessor of the result given in Theorem \ref{Lar2013}. The result of Theorem \ref{Lar2013} is by essentially a logarithmic factor better than the result in \cite{HofLar}.
The reason for this is, that in \cite{HofLar} and \cite{Kritzer} it is worked with techniques analogous to the methods used by Schmidt in \cite{Schmidt1}, whereas in \cite{Lar2013} it is worked with the more powerful techniques of Beck.
Hence it could be conjectured that by suitably adapting the methods of Beck and of \cite{Lar2013} to the discrete case it should be possible to improve also the result of Kritzer by almost a logorithmic factor. Note that in \cite{KriLeoPill} a component by component construction of such point sets was given.

\begin{op}
\textit{Improve the result of Kritzer cited in Theorem \ref{Kritzer1} above on the discrepancy of Hammersley-good lattice point-sets by almost a logarithmic factor.}
\end{op}

The result of Theorem \ref{Lar2013} probably should be essentially the best possible metrical estimate for the discrepancy of Halton-Kronecker sequences.
However, whereas Beck was able to prove this assumption for pure Kronecker sequences, it seems to be out of reach at the moment also to give a rather sharp metrical lower bound in the general case of Halton-Kronecker sequences. The reason for this is that until now we do not even have a satisfying lower bound for the discrepancy of the pure $s$-dimensional Halton sequence in dimension $s \ge 2$ (see Section \ref{discest}).
So (even if it is tempting) we do not state the search for a metrical lower bound for the Halton-Kronecker sequence as an open problem, but just for two special cases

\begin{op}
\textit{Show for the sequence of Problem \ref{vanCorpseq}, namely
$$ (x_n, \{n \sqrt{2}\})_{n \ge 0}, $$
where $(x_n)_{n \ge 1}$ is the van der Corput-sequence in base 2, that 
$$ D_N \ge c \frac{(\log N)^2}{N} $$
holds for a constant $c > 0$ and infinitely many $N$.}
\end{op}

\begin{op}
\textit{Show for the discrepancy $D_N$ of the sequences $(x_n, \{n \alpha\}_{n \ge 0})$ where $(x_n)_{n \ge 0}$ is the van der Corput sequence in base 2, that for almost all $\alpha$ there is a $C(\alpha) > 0$ such that 
$$ D_N \ge C(\alpha) \frac{(\log N)^2}{N} $$
holds for infinitely many $N$.}
\end{op}

The next more general step would be to study the discrepancy of Niederreiter-Kronecker sequences, i.e., hybrid sequences generated be the combination of digital sequences with Kronecker sequences.
Again, like in Problems \ref{copbas} and \ref{copbascomp} it seems to be easier to attack first the problem for \textbf{finite-row}-digital sequences for the digital component.
For the case of infinite-row digital sequences we think it should be challenge enough first to investigate the most basic case. So concerning the analysis of metrical discrepancy of Niederreiter-Kronecker sequences we state the following two open problems:

\begin{op}\label{shupmet}
\textit{Give an essentially sharp upper metrical bound for the discrepancy of Niederreiter-Kronecker sequences
$$ \bsz_n = (\bsx_n, \{n \bsalpha\})_{n \ge 0} $$
where $(\bsx_n)_{n \ge 0}$ is a given digital $s$-dimensional sequence generated by matrices with finite rows, and $(\{n \bsalpha\})_{n \ge 0}$ is a $t$-dimensional Kronecker sequence, which is valid for almost all $\bsalpha$.}
\end{op}

\begin{op}\label{shmetup}
\textit{Give sharp metrical upper - and, if possible, also lower - bounds for the discrepancy of the sequence
$$ \bsz_n = (x_n, \{n \alpha\})_{n \ge 0} $$
where $(\{n \alpha\})_{n \ge 0}$ is a one-dimensional Kronecker sequence, and where $(x_n)_{n \ge 0}$ is the one-dimensional digital sequence in base 2 generated by the matrix
\[
C=
\begin{pmatrix}
1 & 1 & 1 & 1 & 1 & 1 & \ldots\\
0 & 1 & 0 & 0 & 0 & 0 & \ldots\\
0 & 0 & 1 & 0 & 0 & 0 & \ldots\\
0 & 0 & 0 & 1 & 0 & 0 & \ldots\\
\vdots &\vdots & \vdots & \vdots & \vdots & \vdots 
\end{pmatrix}
\in \ZZ_2^{\NN \times \NN} .
\]
The estimates should hold for almost all $\alpha \in \RR$.}
\end{op}

The investigation of the discrepancy of the sequence considered in Problem \ref{shmetup} in a first step leads necessarily to the investigation of the discrepancy of the following one-dimensional sequence:

\begin{op}\label{exmetor}
\textit{Give the exact metrical order of the discrepancy of the sequence
$$ (\{n_k \alpha\})_{k \ge 0} , $$
where $(n_k)_{k \ge 0}$ denotes the increasing sequence of positive integers $n_k$ for which $S_2 (n_k) \equiv 0 \mod{2}$ holds, where $S_2 (n)$ denotes the sum of digits of $n$ in base 2.}
\end{op}

We will make in the remaining part of this section a side-step from discrepancy theory to applications of QMC methods, and especially of hybrid sequences.
The analysis of hybrid sequences already has been motivated by Spanier in \cite{Spanier} and by Keller in \cite{Keller} by applications on transport problems and in image processing. Here we give a suggestion for an application in finance, especially in credit risk management, where hybrid sequences in a quite natural way seemingly should be the suitable tool for generating a simulation scenario. Here by a hybrid sequence, in contrast to the examples given above, we mean the combination of any of the low-discrepancy sequences with a pseudo-random sequence (for the analysis of the distribution of such sequences see for example \cite{GomHofNied}, \cite{Nied2011} or \cite{NiedWin}).

The credit risk management system Credit Metrics of J.P.Morgan (for all theoretical details see \cite{CredMet}) analyses the risk inherent in a large portfolio of credits held by a bank. This is done by calculating the $1 \%$ percentile of the future value of the credit portfolio in one year. To determine this percentile in the original version of credit metrics Monte Carlo simulation is suggested. In this simulation problem each of the credits in the portfolio represents one dimension. Credit portfolios usually contain much more than 1000 single credits, so we have to deal here with a very high-dimensional simulation problem, hence it is certainly too high-dimensional for a pure QMC method. (The benefits of QMC methods in high dimensions usually appear only when the number of sample points is very high.)

On the other hand it could be a good idea to handle the (few) credits in the portfolios which contribute most to the risk of the whole portfolio more carefully, by using a low-discrepancy sequence for the dimensions in the simulation problem which correspond to these highest risk credits. The risk contribution of a credit is a function of the face amount of the credit, its rating class (i.e., of its downgrading-or even default-probabilities), and also on its correlation properties within the credit portfolio. So the (sketch of a) program would be the following: Given a number $N$ of sample points and the dimension $d$ of the simulation problem (i.e., $d$ is the number of credits in the portfolio), 
\begin{itemize}
\item[-] determine a dimension $s < d$, such that it is preferable to use an $s$-dimensional low-discrepancy sequence consisting of $N$ points for an $s$-dimensional simulation problem instead of a pseudo-random sequence,
\item[-] determine the $s$ credits of the portfolio which contribute most to the risk of the whole credit portfolio,
\item[-] carry out the simulations suggested by the system Credit Metrics, but using for the $s$ selected credits a QMC sequence and for all other (many) credits a pseudo-random sequence, i.e., choose a hybrid sequence for the simulation.
\end{itemize}

\begin{op}
\textit{Carry out the above sketched program for using hybrid sequences in the credit risk management program Credit Metrics in detail and analyse the performance in comparison with pure Monte Carlo or pure quasi-Monte Carlo methods.}
\end{op}

\section{Miscellaneous problems}\label{Miscel}

Of course there exist the big open problems in the theory of uniform distribution and essentially in discrepancy theory, like the most prominent one:

\textit{
\begin{itemize}
\item determine the correct order of the best general lower bound for discrepancy, holding for all sequences in $[0,1)^d$ (with the most important contributions of Beck \cite{Beck2}, Bilyk and Lacey \cite{Bilyk}, \cite{BilykL} or Roth \cite{Roth}).\end{itemize}
Or the still open question:
\begin{itemize}
\item is the sequence $\left(\left\{\left(\frac{3}{2}\right)^n \right\}\right)_{n \ge 0}$ uniformly distributed in $[0,1)$ or not?\end{itemize}
Another well-known example is:
\begin{itemize}
\item do there exist $\alpha, \beta \in \RR$ such that the $2$-dimensional Kronecker sequence
$$ (\{n \alpha\}, \{n \beta\})_{n \ge 0} $$
is a low-discrepancy sequence, i.e., satisfies
$$ D_N \le C \frac{(\log N)^2}{N} \quad \mbox{for all} \quad N , $$
or not?
\end{itemize}}

This question has narrow connections to the still open conjecture of Littlewood in diophantine approximation, stating that for all $\alpha, \beta \in \RR$ we have
$$ \underset{n \To \infty}{\lim \inf} \quad n \left\| n \alpha \right\| \cdot \left\| n \beta \right\| = 0, $$
where $\left\| x \right\|$ denotes the distance of $x$ to the nearest integer.\\

In the following we give some further problems in discrepancy theory, which in our opinion also are of considerable interest, and, most probably, will be easier to attack than the above mentioned prominent problems:

The first open problem in this section concerns the best possible lower bound for the star-discrepancy $D_N^*$ for one-dimensional sequences in $[0,1)$.
The until now best known lower bound is already quite old: Bejian \cite{Bejian} has shown in 1982 that for every sequence $(\bsx_n)_{n \ge 0}$ in $[0,1)$ we have 
$$D_N^* \ge \left(0.06015\ldots\right) \frac{\log N}{N}$$ 
for infinitely many $N$.

In \cite{Ostromoukhov} it was shown by Ostromoukhov that there exist $(\bsx_n)_{n \ge 0}$ in $[0,1)$ with 
$$D_N^* \le \left(0.222\ldots\right) \frac{\log N}{N}$$
for all $N$ large enough. So we state

\begin{op}
\textit{Let $c$ be maximal such that for every sequence $(\bsx_n)_{n \ge 0}$ in $[0,1)$ we have
$$ D_N^* \ge c \frac{\log N}{N} $$
for infinitely many $N$.\\
\\
Improve the best until now known bounds $0.06015 \ldots \le c \le 0.222 \ldots$ for $c$.}
\end{op}

First new investigations of the author in this direction show, that by refining a method developed by Liardet \cite{Liardet} and by Tijdeman and Wagner \cite{TijdWag} a rather simple proof for $c \ge 0.06182$ should be possible. 

As already noted above it is the most well-known conjecture in the theory of irregularities of distribution, that for every sequence $(\bsz_n)_{n \ge 0}$ in $[0,1)^d$ we have
$$ D_N \ge c_d \frac{(\log N)^d}{N} $$
for a constant $c_d > 0$ and for infinitely many $N$.
At the moment the best known general lower bound for the discrepancy of sequences in $[0,1)^d$ is
$$ D_N \ge c_d \frac{(\log N)^{\frac{d}{2} + \epsilon(d)}}{N} $$
for infinitely many $N$, with some small $\epsilon(d) > 0$.

Indeed the problem is that even for some of the most well known and seemingly simple sequences in dimensions $d \ge 2$ we do not know the right order of discrepancy, like for example for the $2$-dimensional Kronecker sequences, for the $2$-dimensional Halton sequences and for most of the $2$-dimensional digital low-discrepancy sequences.

Faure has shown in \cite{Faure} that for a certain digital low-discrepancy sequence (i.e., with $D_N \le c \cdot \frac{(\log N)^2}{N}$ for all $N$) we indeed also have
$$ D_N \ge c' \frac{(\log N)^2}{N} $$
for infinitely many $N$.

\begin{op}
\textit{Show that for all $2$-dimensional digital sequences in base 2 we have 
$$ D_N \ge c \frac{(\log N)^2}{N}  $$
for infinitely many $N$.\\
\\
Or, show this estimate at least for a certain class of such (low-discrepancy) sequences (e.g. for NUT sequences or finite-row sequences).}
\end{op}

A probably even more difficult and challenging but for us utmost appealing problem is to find the right order or discrepancy for the Halton sequence in dimension 2. Until now for the discrepancy of this sequence there is no better lower bound known than the general lower bound given by Bilyk and Lacey which holds for \textbf{all} sequences in $[0,1)^2$.

So we finish this paper which should give the reader a survey on some recent results on discrepancy theory and on some open problems in current research with the author's favorite open problem listed in this paper:

\begin{op}
\textit{Give an improved lower bound for the discrepancy of the Halton sequence in dimension 2.\\
\\
In the best case decide whether the right order of the discrepancy is $\frac{(\log N)^2}{N}$, or not.}
\end{op}


\begin{thebibliography}{10}
\bibitem{Beck2} Beck, J.: A two-dimensional van Aardenne Ehrenfest theorem in irregularities of distribution. Compos. Math. (72), pp. 269 - 339. (1989)
\bibitem{Beck} Beck, J.: Probabilistic diophantine approximation, I. Kronecker sequences. Annals of Mathematics (140), pp. 451 - 502. (1994)
\bibitem{Bejian} Bejian, R.: Minoration de la discrepance d'une suite quelconque sur T . Acta Arith. (41), pp. 185 - 202. (1982)
\bibitem{CredMet} Bhatia, M., Finger, Ch.C. and Gupton, G.M.: CreditMetrics - Technical Document. J.P. Morgan, New York. (1997)
\bibitem{Bilyk} Bilyk, D. and Lacey, M.T.: On the small ball inequality in three dimensions. Duke Math. J. (143), pp. 81 -115. (2008)
\bibitem{BilykL} Bilyk, D. and Lacey, M.T.: The supremum norm of the discrepancy function: recent results and connections. To appear in Proc. MCQMC 2012.
\bibitem{Bilyk et} Bilyk, D., Lacey, M.T., Vagharshakyan, A.: On the small ball inequality in all dimensions. J. Funct. Anal. (254), pp. 2470 - 2502. (2008)
\bibitem{BouKont} Bourgain, J. and Kontorovich, A.: On Zaremba's conjecture. arXiv:1107.3776 (2011)
\bibitem{Bykov} Bykovskii,V. A.: The discrepancy of Korobov lattice points. Izvestiya: Mathematics (76), Nr. 3, pp. 446 - 465. (2012)
\bibitem{Dick} Dick, J. and Pillichshammer, F.: Digital sequences, discrepancy theory and quasi-Monte Carlo integration. Cambridge University Press, Cambridge. (2010)
\bibitem{Faure} Faure, H.: Discrepancy lower bound in two dimensions. In: Monte Carlo and Quasi-Monte Carlo Methods in Scientific Computing. Lecture Notes in Statist. (106),
Springer, pp. 198 - 204. (1995)
\bibitem{Frolenkov} Frolenkov, D.A. and Kan, I.D.: A reinforcement of the Bourgain-Kontorovich's theorem by elementary methods II. arXiv:1303.3968v2. (2013)
\bibitem{GomHofNied} Gomez-Perez, D., Hofer, R. and Niederreiter, H.: A general discrepancy bound for hybrid sequences involving Halton sequences. Unif. Distrib. Theory, 8, pp. 31 - 45. (2013)
\bibitem{Hell} Hellekalek, P.: A general discrepancy estimate based on $p$-adic arithmetics. Acta Arith. (139), pp. 117 - 129. (2009)
\bibitem{Hell2} Hellekalek, P.: A notion of diaphony based on $p$-adic arithmetic. Acta Arith. (145), pp. 273 - 284. (2010)
\bibitem{HellNied} Hellekalek, P. and  Niederreiter, H.: Constructions of uniformly distributed sequences using the $b$-adic method. Uniform Distribution Theory (6), pp. 185 - 200. (2011)
\bibitem{Hofer} Hofer, R.: On the distribution of Niederreiter-Halton sequences. J. Number Theory, (129), pp. 451-463. (2009) 
\bibitem{HofKriLarPill} Hofer, R., Kritzer, P., Larcher G. and Pillichshammer, F.: Distribution properties of generalized van der Corput-Halton sequences and their subsequences. Int. J. Number Theory (5), pp. 719 - 746. (2009)
\bibitem{HofLar} Hofer, R. and Larcher, G.: On existence and discrepancy of certain digital Niederreiter-Halton sequences. Acta Arith. (141), pp. 369 - 394. (2010)
\bibitem{HofLar2} Hofer, R. and Larcher, G.: Metrical results on the discrepancy of Halton-Kronecker sequences. Math. Zeit. (271), pp. 1 - 11. (2012)
\bibitem{Keller} Keller, A.: Quasi-Monte Carlo image synthesis in a nutshell. Submitted. (2012)
\bibitem{Kontorovich} Kontorovich, A.: From Apollonius to Zaremba: Local-global phenomena in thin orbits. Bull. Amer. Math. Soc. (50), pp.187 - 228. (2013)
\bibitem{Kritzer} Kritzer, P.: On an example of finite mixed quasi-Monte Carlo point sets. Monatsh. Math. (168), pp. 443 - 459. (2012)
\bibitem{KriLeoPill} Kritzer, P., Leobacher, G. and Pillichshammer, F.: Component-by-component construction of hybrid point sets based on Hammersley and lattice point sets. Monte Carlo and Quasi-Monte Carlo Methods 2012, to appear. (2013)
\bibitem{KriPill} Kritzer, P. and Pillichshammer, F.: Low discrepancy polynomial lattice point sets. J. Number Theory (132), pp. 2510 - 2534. (2012)
\bibitem{KuNied2006} Kuipers, L. and Niederreiter, H.: {\it{Uniform distribution of sequences}}. John Wiley, New York, 1974. Reprint, Dover Publications, Mineola, NY. (2006) 
\bibitem{Lar1986} Larcher, G.: On the distribution of sequences connected with good lattice points, Monatsh. Math. 101, pp. 135 - 150. (1986)
\bibitem{Lar1993} Larcher, G.: Nets obtained from rational functions over finite fields. Acta Arith. (63), pp. 1 - 13. (1993)
\bibitem{Lar1995} Larcher, G.: On the distribution of an analog to classical Kronecker-sequences. J. Number Theory 52, 198 - 215. (1995)
\bibitem{Lar1998} Larcher, G.: Digital point sets: analysis and application. In: Random and quasi-random point sets. vol. 138 of Lecture Notes in Statistics, pp. 167 - 222. Springer, New York (1998)
\bibitem{Lar-Pill} Larcher, G.: On the distribution of digital sequences. In: Monte Carlo and Quasi-Monte Carlo Methods 1996 (Salzburg), vol.127 of Lecture Notes in Statistics, pp. 109 - 123, Springer, New York. (1998)
\bibitem{Lar2013} Larcher, G.: Probabilistic Diophantine approximation and the distribution of Halton-Kronecker sequences. J. Complexity, to appear. (2013)
\bibitem{Lar2} Larcher, G. and Niederreiter, H.: Kronecker-type sequences and nonarchimedean diophantine approximation. Acta Arith. 63, 380 - 396. (1993)
\bibitem{Lar} Larcher, G. and Niederreiter, H.: Generalized $(t,s)$-sequences, Kronecker-type sequences, and diophantine approximation of formal Laurent series. Trans. Amer. Math. Soc. (347), pp. 2051 - 2073. (1995)
\bibitem{Lar-Pill2} Larcher, G. and Pillichshammer, F.: A metrical best possible lower bound on the star discrepancy of digital sequences. Monatsh. Math., to appear. (2013)
\bibitem{Lar-Pill3} Larcher, G. and Pillichshammer, F.: Metrical lower bounds on the discrepancy of digital Kronecker sequences. Submitted for publication. arXiv:1302.5267.  (2013)
\bibitem{Levin} Levin, M.B.: On low discrepancy sequences and low discrepancy ergodic transformations of the multidimensional unit cube. Israel Journal of Mathematics (178), pp. 61 - 106. (2010)
\bibitem{Liardet} Liardet, P.: Discrepance sur le cercle. Primaths I, Univ. Marseille , pp. 7 - 11. (1979)
\bibitem{Nied1} Niederreiter, H.: Point sets and sequences with small discrepancy. Monatsh.Math. (104), pp. 273 - 337. (1987)
\bibitem{Nied2} Niederreiter, H.: Random number generation and quasi-Monte Carlo methods. No. 63 in CBMS-NSF Series in Applied Mathematics, SIAM, Philadelphia. (1992)
\bibitem{Nied2009} Niederreiter, H.: On the discrepancy of some hybrid sequences. Acta Arith. (138), pp. 373 - 398. (2009)
\bibitem{Nied2010} Niederreiter, H.: Further discrepancy bounds and an Erd\"os-Turan-Koksma inequality for hybrid sequences. Monatsh. Math. (161), pp. 193 - 222. (2010)
\bibitem{Nied2011} Niederreiter, H.: Discrepancy bounds for hybrid sequences involving matrix-method pseudorandom vectors. Publ. Math. Debrecen, 79, pp. 589 - 603. (2011)
\bibitem{Nied2012} Niederreiter, H.: Improved discrepancy bounds for hybrid sequences involving Halton sequences. Acta Arith. (155), pp. 71 - 84. (2012)
\bibitem{NiedWin} Niederreiter, H. and Winterhof, A.: Discrepancy bounds for hybrid sequences involving digital explicit inversive pseudorandom numbers. Unif. Distrib. Theory, 6, pp. 33 - 56. (2011)
\bibitem{Ostromoukhov} Ostromoukhov, V.: Recent progress in improvement of extreme discrepancy and star discrepancy of one-dimensional sequences. In: Monte-Carlo and Quasi-Monte Carlo Methods 2008, P. L'Ecuyer and A. B. Owen (eds.), Springer, New York, pp. 561 - 572. (2009)
\bibitem{Roth} Roth, K. F.: On irregularities of distribution. Mathematika 1, pp. 73 - 79. (1954)
\bibitem{Schmidt} Schmidt, W. M.: A metrical theorem in diophantine approximation. Canad. J. Math. (12), pp. 619 - 631. (1960)
\bibitem{Schmidt1} Schmidt, W.M.: Metrical theorems on fractional parts of sequences. Trans. Amer. Math. Soc. (110), pp. 493 - 512. (1964)
\bibitem{Schmidt2} Schmidt, W.M.: Irregularities of distribution VII. Acta Arith. (66), 21, pp. 45 - 50. (1972)
\bibitem{Spanier} Spanier, J.: Quasi-Monte Carlo methods for particle transport problems, in: H. Niederreiter, P.J.-S. Shiue (Eds.), Monte Carlo and Quasi-Monte Carlo Methods in Scientific Computing, in: Lecture Notes in Statistics, vl.106, Springer, New York, pp. 121 - 148. (1995)
\bibitem{TijdWag} Tijdeman, R. and Wagner, G.: A sequence has almost nowhere small discrepancy. Monatsh. Math. (90), pp. 315 - 329. (1980)
\end{thebibliography}
\end{document}